# DATA-DRIVEN SOBOLEV TESTS OF UNIFORMITY ON COMPACT RIEMANNIAN MANIFOLDS

By P. E. Jupp

*University of St. Andrews*

Data-driven versions of Sobolev tests of uniformity on compact Riemannian manifolds are proposed. These tests are invariant under isometries and are consistent against all alternatives. The large-sample asymptotic null distributions are given.

**1. Introduction.** A fundamental hypothesis in directional statistics is that of uniformity of a distribution on a compact Riemannian manifold, such as the circle, sphere or rotation group. An important large class of tests of uniformity consists of Giné's [7] Sobolev tests. Each Sobolev test is specified by a sequence $a_1, a_2, \ldots$ of real numbers satisfying a suitable square-summability condition, which is given in (2.1) below. Sobolev tests for which only a few $a_k$ are nonzero are often simple to calculate, whereas those for which all $a_k$ are nonzero are consistent against all alternatives. The vastness of the class of Sobolev tests presents a problem, as the statistician needs to choose the sequence $a_1, a_2, \ldots$. This paper solves the problem by providing a simple data-driven version of Sobolev tests of uniformity on arbitrary compact connected Riemannian manifolds, in which the sequence $a_1, a_2, \ldots$ is not chosen in advance, but is specified automatically by the data. These data-driven tests are invariant under isometries and consistent against all alternatives.

The construction of the data-driven tests is based on the fact that Sobolev tests of uniformity can be regarded as weighted score tests of uniformity within some of the canonical exponential models introduced by Beran [1]. Thus, they are analogous to Neyman's [21] smooth tests of uniformity on the unit interval. Data-driven versions of Neyman's smooth tests were introduced by Ledwina [18] and refined in Inglot, Kallenberg and Ledwina [12] and Kallenberg and Ledwina [16, 17]. Consistency and other desirable properties of these tests were established by Kallenberg and Ledwina [15], Inglot, Kallenberg and Ledwina [12], Inglot [11] and Inglot and Ledwina [13, 14]. By suitable adaptation of the version of Ledwina's test considered in [12, 16] and [17], Bogdan, Bogdan and Futschik [2] obtained a data-driven test of uniformity on the circle. The data-driven tests presented in this paper generalize those of Bogdan et al. from the circle to arbitrary compact Riemannian manifolds.

In Section 2, the construction of Sobolev tests of uniformity is reviewed and the data-driven versions are introduced. Results on consistency and other large-sample asymptotic properties are given in Section 3. Section 4 gives details in the important cases in which the sample space is a sphere, a projective space or the rotation group $SO(3)$. Some simulation results are reported in Section 5.

An outline of the requisite differential geometry can be found in Section 2 of [7] or the Appendix of [9]. More detailed accounts are given in [8] and [3].

## 2. Data-driven Sobolev tests of uniformity.

2.1. *Sobolev tests of uniformity.* Let $M$ be a compact connected Riemannian manifold without boundary. The Riemannian metric determines the uniform probability measure $\mu$ on $M$. The intuitive idea behind Giné's [7] Sobolev tests of uniformity is to map the manifold $M$ into the Hilbert space $L^2(M, \mu)$ of square-integrable real-valued functions on $M$ by a function $\mathbf{t}: M \to L^2(M, \mu)$ such that if $x$ is uniformly distributed then $E[\mathbf{t}(x)] = \mathbf{0}$, and to reject uniformity if the sample mean of $\mathbf{t}(x)$ is "far" from $\mathbf{0}$.

The standard way of constructing such mappings $\mathbf{t}$ is due to Giné [7] and is based on the eigenfunctions of the Laplacian operator on $M$ acting on smooth real-valued functions on $M$. For $k \geq 1$, let $E_k$ denote the space of eigenfunctions corresponding to $\lambda_k$, the $k$th nonzero eigenvalue, and put $d_k = \dim E_k$. Then there is a well-defined map, $\mathbf{t}_k$ of $M$ into $E_k$ given by

$$\mathbf{t}_k(x) = \sum_{i=1}^{d_k} f_i(x) f_i,$$

where $\{f_i : 1 \leq i \leq d_k\}$ is any orthonormal basis of $E_k$. If $a_1, a_2, \ldots$ is a sequence of real numbers such that

$$\sum_{k=1}^{\infty} a_k^2 d_k < \infty \tag{2.1}$$

then

$$x \mapsto \mathbf{t}(x) = \sum_{k=1}^{\infty} a_k \mathbf{t}_k(x) \tag{2.2}$$

defines a mapping $\mathbf{t}$ of $M$ into $L^2(M,\mu)$. The resulting Sobolev statistic evaluated on observations $x_1,\ldots,x_n$ on $M$ is

$$T_n = \frac{1}{n}\left\|\sum_{i=1}^n \mathbf{t}(x_i)\right\|^2 = \frac{1}{n}\sum_{i=1}^n\sum_{j=1}^n \langle \mathbf{t}(x_i), \mathbf{t}(x_j)\rangle,$$

where $\langle,\rangle$ denotes the inner product on $L^2(M,\mu)$ given by

$$\langle f, g\rangle = \int_M f(x)g(x)\,d\mu(x).$$

The corresponding Sobolev test rejects uniformity for large values of $T_n$.

2.2. *Score tests of uniformity.* Any vector $(\boldsymbol{\theta}_1,\ldots,\boldsymbol{\theta}_k)$ in $\bigoplus_{j=1}^k E_j$ determines a distribution on $M$ with density (with respect to the uniform distribution)

$$(2.3) \qquad f(\mathbf{x};\boldsymbol{\theta}_1,\ldots,\boldsymbol{\theta}_k) = \exp\left\{\sum_{j=1}^k \langle \boldsymbol{\theta}_j, \mathbf{t}_j(\mathbf{x})\rangle - \kappa(\boldsymbol{\theta}_1,\ldots,\boldsymbol{\theta}_k)\right\},$$

where $\kappa(\boldsymbol{\theta}_1,\ldots,\boldsymbol{\theta}_k)$ is a log normalizing constant. The class of such distributions was introduced by Beran [1]. The quadratic score statistic $S_k$ in the score test of uniformity (i.e., $\boldsymbol{\theta}_j = \mathbf{0}$ for $j = 1,\ldots,k$) in the exponential model (2.3) is

$$(2.4) \qquad S_k = n\left\|\frac{1}{n}\sum_{i=1}^n \mathbf{t}_{(k)}(x_i)\right\|^2,$$

where $\mathbf{t}_{(k)}$ is $\mathbf{t}$ defined by (2.2) with

$$(2.5) \qquad a_j = \begin{cases} 1, & \text{for } j \leq k, \\ 0, & \text{for } j > k. \end{cases}$$

Thus $S_k$ is the Sobolev statistic corresponding to (2.5). These $\mathbf{t}_{(k)}$ are used in Hendriks's [9] density estimates $\hat{f}_k$, which are given by

$$\hat{f}_k(x) = \frac{1}{n}\sum_{i=1}^n \langle \mathbf{t}_{(k)}(x_i), \mathbf{t}_{(k)}(x)\rangle.$$

2.3. *Data-driven tests of uniformity.* A major problem in practice with the score tests of Section 2.2 is the need to choose a suitable $k$. The solution proposed here is for this choice to be made by the data, using a modification of Schwarz's [24] Bayesian Information Criterion (BIC) selection rule, as considered by Inglot et al. [12], Kallenberg and Ledwina [16, 17], Inglot [11] and Inglot and Ledwina [14] for observations on the line and by Bogdan et al. [2] for observations on the circle.

The choice of $k$ in the data-driven versions of the score tests of Section 2.2 is based on the *penalized score statistic*

$$(2.6) \qquad B_S(k) = S_k - \nu_k \log n,$$

where

$$\nu_k = \sum_{i=1}^{k} d_i.$$

The second term on the right-hand side of (2.6) penalizes higher-dimensional models. The value of $k$ is chosen as $\hat{k}$, where

$$(2.7) \qquad \hat{k} = \inf\left\{k \in \mathbb{N} : B_S(k) = \sup_{m \in \mathbb{N}} B_S(m)\right\}.$$

(Recall that $\inf \varnothing = \infty$.) This procedure is an analogue for score tests of Schwarz's [24] BIC selection rule. The data-driven score tests reject the null hypothesis of uniformity for large values of $S_k$.

The main large-sample properties of these tests can be stated intuitively as:

(i) Under the null hypothesis, $\hat{k}$ tends to be near 1, so that the test statistic $S_k$ is "often" simple (for large samples);

(ii) The asymptotic null distribution of $S_k$ is chi-squared;

(iii) Under the alternative hypothesis, $\hat{k}$ tends to rise to a value large enough to cause rejection of the null hypothesis.

Rigorous versions of (i), (ii) and (iii) are given in Theorems 3.1, 3.2 and 3.3, respectively.

A useful geometrical tool in deriving the properties of data-driven score tests is the *spectral function*, which can be defined as

$$(2.8) \qquad e(x, y, T) = \frac{1}{\mathrm{vol}(M)} \sum_{\lambda_k \leq T} \langle \mathbf{t}_k(x), \mathbf{t}_k(y) \rangle, \qquad x, y \in M, T > 0,$$

where $\mathrm{vol}(M)$ denotes the Riemannian volume of $M$. Note that

$$(2.9) \qquad e(x, x, \lambda_k) = \frac{1}{\mathrm{vol}(M)} \|\mathbf{t}_{(k)}(x)\|^2.$$

A key property of the spectral function is

$$(2.10) \qquad \sup_{x \in M} |T^{-m/2} e(x, x, T)(2\sqrt{\pi})^m \Gamma(m/2 + 1) - 1| = O(T^{-1/2}),$$

$$T \to \infty,$$

where $m$ is the dimension of $M$ and $\Gamma$ denotes the gamma function. See, for example Theorem 1.1 of [10] or equation (2.25) of [6]. Let $N(T)$ be the number of eigenvalues $\lambda$ (counted with their multiplicities) of the Laplacian with $\lambda \leq T$. Then

$$\int_M e(x, x, T) \, \mathrm{dvol}(x) = \mathrm{vol}(M) N(T),$$

and so (2.10) yields Weyl's formula

(2.11) $$\lim_{T \to \infty} T^{-m/2} N(T) = \frac{\mathrm{vol}(M)}{(2\sqrt{\pi})^m \Gamma(m/2 + 1)}.$$

See, for example, page 243 of [20] or page 9 of [4]. Thus

(2.12) $$\nu_k \sim \frac{\mathrm{vol}(M)}{(2\sqrt{\pi})^m \Gamma(m/2 + 1)} \lambda_k^{m/2} \qquad \text{as } k \to \infty.$$

Combining (2.10) with (2.11) gives

(2.13) $$\frac{\mathrm{vol}(M) e(x, x, \lambda_k)}{\nu_k} \to 1 \qquad \text{(uniformly in } x\text{) as } k \to \infty.$$

The following complement of (2.13) will be useful.

PROPOSITION 2.1. *If $x \neq y$ then*

(2.14) $$\frac{e(x, y, \lambda_k)}{\nu_k} \to 0 \qquad \text{as } k \to \infty.$$

Since the proof of Proposition 2.1 is rather technical, it is postponed to the Appendix.

If $\hat{k} = \infty$, then $S_k$ is not defined. The following proposition shows that this occurs with probability zero, except in very small samples.

PROPOSITION 2.2. *If $n \geq 3$ then for random samples of size $n$ from a continuous distribution on $M$,*

(2.15) $$P(\hat{k} = \infty) = 0.$$

PROOF. Since the distribution is continuous, the observations $x_1, \ldots, x_n$ are distinct with probability 1. In this case, (2.4), (2.8), (2.9), (2.13), and (2.14) give

$$\begin{aligned}
\frac{S_k}{\nu_k} &= \frac{1}{n\nu_k} \left\{ \sum_{i=1}^n \|\mathbf{t}_{(k)}(x_i)\|^2 + \sum_{i=1}^n \sum_{j \neq i} \langle \mathbf{t}_{(k)}(x_i), \mathbf{t}_{(k)}(x_j) \rangle \right\} \\
&= \frac{\mathrm{vol}(M)}{n} \left\{ \sum_{i=1}^n \frac{e(x_i, x_i, \lambda_k)}{\nu_k} + \sum_{i=1}^n \sum_{j \neq i} \frac{e(x_i, x_j, \lambda_k)}{\nu_k} \right\} \\
&\to 1 \qquad \text{as } k \to \infty.
\end{aligned}$$

It follows that

$$\frac{B_S(k) - B_S(1)}{\nu_k} \to 1 - \log n \quad \text{as } k \to \infty.$$

Thus, if $n \geq 3$ then $B_S(k) < B_S(1)$ for large enough $k$, so that $\hat{k} < \infty$. □

REMARK. In practice $\hat{k}$ is calculated not by (2.7) but as

$$\inf\left\{k : B_S(k) = \sup_{1 \leq m \leq K} B_S(m)\right\}$$

for some suitable $K$. Tables 1, 3 and 5 in Section 5 indicate that it is reasonable to take $K = 5$ for $M = S^2, \mathbb{R}P^2$ or $SO(3)$. Calculation of $S_{\hat{k}}$ takes approximately $K$ times as much effort as calculation of a Sobolev test which is consistent against all alternatives. Thus, the extra computational cost of using a data-driven test is small.

## 3. Asymptotic properties.

THEOREM 3.1. *Under uniformity,*

$$\lim_{n \to \infty} P(\hat{k} = 1) = 1.$$

PROOF. Since $\{f_i : 1 \leq i \leq d_k\}$ is a orthonormal basis of $E_k$, if $x$ is uniformly distributed then

(3.1) $$E[\mathbf{t}_{(k)}(x)] = \mathbf{0} \quad \text{and} \quad E[\mathbf{t}_{(k)}(x)\mathbf{t}_{(k)}(x)'] = \mathbf{I}_{\nu_k}.$$

Straight forward calculation gives

(3.2) $$E[S_k] = \nu_k$$

and

$$\text{var}(S_k) = \frac{1}{n^2}\text{var}\left(\sum_{i=1}^n \|\mathbf{t}_{(k)}(x_i)\|^2 + \sum_{i=1}^n \sum_{j \neq i} \langle \mathbf{t}_{(k)}(x_i), \mathbf{t}_{(k)}(x_j) \rangle\right)$$

$$= \frac{1}{n^2} \sum_{i=1}^n \text{var}(\|\mathbf{t}_{(k)}(x_i)\|^2)$$

$$+ \frac{2}{n^2} \sum_{i=1}^n \sum_{j \neq i} \text{tr}(E[\mathbf{t}_{(k)}(x_i)\mathbf{t}_{(k)}(x_i)']E[\mathbf{t}_{(k)}(x_j)\mathbf{t}_{(k)}(x_j)'])$$

(3.3) $$= \frac{1}{n}\text{var}(\|\mathbf{t}_{(k)}(x)\|^2) + 2\left(1 - \frac{1}{n}\right)\nu_k.$$

It follows from (2.9)–(2.12) that there is a positive $A$ such that

$$|\|\mathbf{t}_{(k)}(x)\|^2 - \nu_k| \leq A\nu_k^{(m-1)/m} \qquad \text{for } k = 1, 2, \ldots \text{ and } x \in M.$$

Thus

$$\operatorname{var}(\|\mathbf{t}_{(k)}(x)\|^2) \leq A^2 \nu_k^{2(1-1/m)} \qquad \text{for } k = 1, 2, \ldots.$$

Combining this with (3.3) gives

$$\operatorname{var}(S_k) < A^2 \nu_k^{2(1-1/m)} + 2\nu_k. \tag{3.4}$$

It follows from (2.11) that there is a positive $T_0$ such that

$$T > T_0 \implies 2 \leq N(3^{2/m}T)/N(T) \leq 4. \tag{3.5}$$

Put $\mu_1 = \lambda_1$, choose $\mu_2$ such that $N(\mu_2) = \nu_2$, and define $\mu_3, \mu_4, \ldots$ by

$$\mu_l = 3^{2(l-1)/m} \mu_2 \qquad \text{for } l \geq 3. \tag{3.6}$$

From (3.5) and (3.6) there are positive constants $C_1$ and $C_2$ such that

$$N(\mu_l) \geq 2^l C_1 \qquad \text{for } l \geq 1, \tag{3.7}$$

$$N(\mu_l)/N(\mu_{l-1}) \leq C_2 \qquad \text{for } l \geq 2. \tag{3.8}$$

Suppose that $\mu_l \leq \lambda_{\hat{k}} < \mu_{l+1}$, for some $l \geq 2$. Define $r$ by $N(\mu_{l+1}) = \nu_r$. Since $k \mapsto S_k$ and $N$ are increasing functions,

$$S_r - N(\mu_l) \log n \geq S_{\hat{k}} - \nu_{\hat{k}} \log n > S_1 - \nu_1 \log n,$$

so that, by (3.2),

$$\begin{aligned}
S_r - E[S_r] &\geq S_1 - \nu_1 \log n + N(\mu_l) \log n - N(\mu_{l+1}) \\
&\geq (N(\mu_l) - \nu_1) \log n - N(\mu_{l+1}) \\
&= N(\mu_l)\{[1 - \nu_1/N(\mu_l)] \log n - N(\mu_{l+1})/N(\mu_l)\} \\
&\geq N(\mu_l)\{[d_2/(d_1 + d_2)] \log n - C_2\},
\end{aligned} \tag{3.9}$$

which is positive if $\log n > C_2(d_1 + d_2)/d_2$. In this case, it follows from (3.9), Chebyshev's inequality and (3.4) that

$$\begin{aligned}
&P(\mu_l \leq \lambda_{\hat{k}} < \mu_{l+1}) \\
&\leq \frac{\operatorname{var}(S_r)}{\{N(\mu_l)[d_2/(d_1+d_2) \log n - C_2]\}^2} \\
&< \left\{\frac{A^2}{N(\mu_{l+1})^{2/m}}\left(\frac{N(\mu_{l+1})}{N(\mu_l)}\right)^2 + \frac{2}{N(\mu_l)}\frac{N(\mu_{l+1})}{N(\mu_l)}\right\}\left(\frac{d_2}{d_1+d_2}\log n - C_2\right)^{-2}.
\end{aligned}$$

Provided that $\log n > 2C_2(d_1 + d_2)/d_2$, using (3.8) then gives

$$(3.10) \quad P(\mu_l \leq \lambda_{\hat{k}} < \mu_{l+1}) < \left\{ \frac{A^2 C_2^2}{N(\mu_{l+1})^{2/m}} + \frac{2C_2}{N(\mu_l)} \right\} \frac{4(d_1 + d_2)^2}{d_2^2 (\log n)^2}.$$

Summing (3.10) over $l = 2, 3, \ldots$ and using (3.7) and (2.15) gives

$$P(\hat{k} > 1) < C_3 (\log n)^{-2}$$

for some positive $C_3$, and so

$$P(\hat{k} > 1) \to 0 \quad \text{as } n \to \infty. \qquad \square$$

The asymptotic null distribution of $S_{\hat{k}}$ can now be found.

THEOREM 3.2. *Under uniformity,*

$$S_{\hat{k}} \xrightarrow{d} \chi^2_{\nu_1} \quad \text{as } n \to \infty,$$

*where $\xrightarrow{d}$ denotes convergence in distribution and $\chi^2_{\nu_1}$ denotes the chi-squared distribution with $\nu_1$ degrees of freedom.*

PROOF. This is a consequence of Theorem 3.1 and the fact [which can be obtained by applying the central limit theorem to $\mathbf{t}_{(k)}(x)$ and using (3.1)] that, under uniformity, for any fixed $k$,

$$S_k \xrightarrow{d} \chi^2_{\nu_k} \quad \text{as } n \to \infty. \qquad \square$$

The following theorem guarantees consistency.

THEOREM 3.3. *The test which rejects uniformity for large values of $S_{\hat{k}}$ is consistent against all alternatives to uniformity.*

PROOF. For any nonuniform distribution, there is a natural number $K$ such that $E[\mathbf{t}_K(x)] \neq \mathbf{0}$. For $1 \leq j < K$,

$$\begin{aligned} P(\hat{k} = j) &\leq P(S_j - \nu_j \log n \geq S_K - \nu_K \log n) \\ &\leq P(n \|\bar{\mathbf{t}}_K\|^2 \leq (\nu_K - \nu_j) \log n) \\ &\to 0 \quad \text{as } n \to \infty, \end{aligned}$$

because $\|\bar{\mathbf{t}}_K\|^2 \to \|E[\mathbf{t}_K(x)]\|^2$ almost surely as $n \to \infty$ and $(\log n)/n \to 0$ as $n \to \infty$. Then, for any positive $C$,

$$\begin{aligned} P(S_{\hat{k}} > C) &\geq P(S_K > C) - P(\hat{k} < K) \\ &\geq P(n \|\bar{\mathbf{t}}_K\|^2 > C) - P(\hat{k} < K) \\ &\to 1 \quad \text{as } n \to \infty. \end{aligned}$$

Since $S_{\hat{k}}$ has a nondegenerate limiting distribution under uniformity, it follows that the test is consistent. $\square$

**4. Examples: spheres, projective spaces and the rotation group.** In order to calculate the statistic $S_{\hat{k}}$, explicit expressions for $\langle \mathbf{t}_k(x), \mathbf{t}_k(y) \rangle$ and $\nu_k$ are required. This section gives such expressions where $M$ is a sphere, a projective space or the rotation group $SO(3)$.

4.1. *Spheres.* It follows from the formula for the cosine of a difference and from Proposition 2.1 of [22] that, for $\mathbf{x}, \mathbf{y}$ in $S^{p-1}$ (the unit sphere in $\mathbb{R}^p$)

$$(4.1) \qquad \langle \mathbf{t}_k(\mathbf{x}), \mathbf{t}_k(\mathbf{y}) \rangle = \begin{cases} 2\cos(k\theta), & \text{if } p = 2, \\ \left(1 + \dfrac{k}{\alpha}\right) C_k^\alpha(\mathbf{x}'\mathbf{y}), & \text{if } p > 2, \end{cases}$$

where $\cos\theta = \mathbf{x}'\mathbf{y}$, $\alpha = p/2 - 1$ and $C_k^\alpha$ denotes the Gegenbauer polynomial of degree $k$. The expression for $d_k$ given on page 171 of [22] yields

$$(4.2) \qquad \nu_k = \frac{1}{p-1}\left\{ k \binom{p+k-2}{p-2} + (k+1)\binom{p+k-1}{p-2} \right\} - 1.$$

In the case of $S^2$, (4.1) reduces to

$$\langle \mathbf{t}_k(\mathbf{x}), \mathbf{t}_k(\mathbf{y}) \rangle = (2k+1) P_k(\mathbf{x}'\mathbf{y})$$

(as in (6.7) of [7]), where $P_k$ denotes the Legendre polynomial of degree $k$, and (4.2) reduces to

$$\nu_k = k(k+2).$$

4.2. *Projective spaces.* Since the eigenspace $E_k$ of the Laplacian on $\mathbb{R}P^{p-1}$ (the projective space of one-dimensional subspaces of $\mathbb{R}^p$), can be identified with the eigenspace $E_{2k}$ of the Laplacian on $S^{p-1}$, it follows from (4.1) that, for $\pm\mathbf{x}, \pm\mathbf{y}$ in $\mathbb{R}P^{p-1}$,

$$(4.3) \qquad \langle \mathbf{t}_k(\pm\mathbf{x}), \mathbf{t}_k(\pm\mathbf{y}) \rangle = \begin{cases} 2\cos(2k\theta), & \text{if } p = 2 \\ \left(1 + \dfrac{2k}{\alpha}\right) C_{2k}^\alpha(\mathbf{x}'\mathbf{y}), & \text{if } p > 2, \end{cases}$$

where $\cos\theta = \mathbf{x}'\mathbf{y}$, and from (4.2) that

$$(4.4) \qquad \nu_k = \sum_{i=1}^{k} \left\{ \binom{p+2i-3}{p-2} + \binom{p+2i-2}{p-2} \right\}.$$

In the case of $\mathbb{R}P^2$, (4.3) and (4.4) reduce to

$$\langle \mathbf{t}_k(\pm\mathbf{x}), \mathbf{t}_k(\pm\mathbf{y}) \rangle = (4k+1) P_{2k}(\mathbf{x}'\mathbf{y})$$

and

$$\nu_k = k(3k+2).$$

4.3. *The rotation group $SO(3)$.* There is a standard identification of the rotation group $SO(3)$ with $\mathbb{R}P^3$ given by the mapping which sends $\pm\mathbf{u} = \pm(u_1,\ldots,u_4)'$ in $\mathbb{R}P^3$ to the matrix

$$\begin{pmatrix} u_1^2 + u_2^2 - u_3^2 - u_4^2 & -2(u_1u_4 - u_2u_3) & 2(u_1u_3 + u_2u_4) \\ 2(u_1u_4 + u_2u_3) & u_1^2 + u_3^2 - u_2^2 - u_4^2 & -2(u_1u_2 - u_3u_4) \\ -2(u_1u_3 - u_2u_4) & 2(u_1u_2 + u_3u_4) & u_1^2 + u_4^2 - u_2^2 - u_3^2 \end{pmatrix}$$

in $SO(3)$. Combining this with (4.3) and (4.4), for $p = 4$ shows that on $SO(3)$

$$\langle \mathbf{t}_k(\mathbf{X}), \mathbf{t}_k(\mathbf{Y}) \rangle = (2k+1)C_{2k}^1\left(\frac{\operatorname{tr}(\mathbf{X}'\mathbf{Y}) + 1}{4}\right), \qquad \mathbf{X}, \mathbf{Y} \in SO(3),$$

$$\nu_k = \frac{k(4k^2 + 12k + 11)}{3}.$$

**5. Simulation study.** Tables 1–6 summarize the null distributions of $\hat{k}$ and $S_{\hat{k}}$ for $S^2$, $\mathbb{R}P^2$ and $SO(3)$, based on 10,000 simulations. Comparison with the simulations for the circle $S^1$ by Bogdan et al. [2] show that, as $n \to \infty$, $\hat{k}$ converges more rapidly to 1 and the distribution of $S_{\hat{k}}$ converges more rapidly to its limiting $\chi^2$ distribution for $S^2$, $\mathbb{R}P^2$ and $SO(3)$ than for $S^1$.

Tables 2, 4 and 6 show that comparison of the observed values of $S_{\hat{k}}$ with the upper 10% and 5% quantiles of the large-sample asymptotic $\chi^2$ distribution is reasonable for $n \geq 30$ for $S^2$, $n \geq 20$ for $\mathbb{R}P^2$ and $n \geq 25$ for $SO(3)$. Since the asymptotic distribution of $S_{\hat{k}}$ is $\chi^2$, it follows from general results of Cordeiro and Ferrari [5] that there are polynomial modifications of $S_{\hat{k}}$ which bring its null distribution closer to the asymptotic distribution. The cubic modification $S_{\hat{k}}^*$ of $S_{\hat{k}}$ given by

$$(5.1) \quad S_{\hat{k}}^* = \left\{1 + \frac{1.37 - 0.31 S_{\hat{k}}}{n}\right\} S_{\hat{k}} \qquad \text{for } M = S^2,$$

$$(5.2) \quad S_{\hat{k}}^* = \left\{1 + \frac{1.91 - 0.21 S_{\hat{k}}}{n}\right\} S_{\hat{k}} \qquad \text{for } M = \mathbb{R}P^2,$$

Table 1
*Empirical distribution of $\hat{k}$ (based on 10,000 simulations) from samples of size n from the uniform distribution on $S^2$*

|   | $n = 5$ | $n = 10$ | $n = 15$ | $n = 20$ | $n = 25$ | $n = 30$ |
|---|---------|----------|----------|----------|----------|----------|
| 1 | 8407 | 9635 | 9826 | 9900 | 9951 | 9964 |
| 2 | 1060 | 338 | 167 | 97 | 49 | 36 |
| 3 | 357 | 25 | 6 | 3 | 0 | 0 |
| 4 | 92 | 2 | 1 | 0 | 0 | 0 |
| 5–10 | 84 | 0 | 0 | 0 | 0 | 0 |

TABLE 2
*Empirical upper tail probabilities $P(S_{\hat{k}} \geq \chi^2_{3;\alpha})$ (upper line) and $P(S^*_{\hat{k}} \geq \chi^2_{3;\alpha})$ (lower line) from samples of size $n$ from the uniform distribution on $S^2$. The estimates are based on 10,000 simulations. The modified statistic $S^*_{\hat{k}}$ is defined in (5.1). Bold figures indicate values in the interval $\alpha \pm 2\sqrt{\alpha(1-\alpha)/10{,}000}$*

| $\alpha$ | $n=5$ | $n=10$ | $n=15$ | $n=20$ | $n=25$ | $n=30$ |
|---|---|---|---|---|---|---|
| 0.10 | 0.232 | 0.132 | 0.114 | 0.108 | 0.110 | **0.102** |
|      | 0.078 | 0.114 | **0.105** | **0.101** | **0.102** | 0.097 |
| 0.05 | 0.187 | 0.078 | 0.062 | 0.058 | 0.058 | **0.051** |
|      | 0.000 | 0.056 | **0.050** | **0.048** | **0.048** | 0.043 |
| 0.01 | 0.120 | 0.043 | 0.023 | 0.017 | 0.014 | **0.011** |
|      | 0.000 | 0.000 | 0.018 | 0.015 | **0.010** | **0.008** |

TABLE 3
*Empirical distribution of $\hat{k}$ (based on 10,000 simulations) from samples of size $n$ from the uniform distribution on $\mathbb{R}P^2$*

|     | $n=5$ | $n=10$ | $n=15$ | $n=20$ | $n=25$ | $n=30$ |
|-----|-------|--------|--------|--------|--------|--------|
| 1   | 9044  | 9900   | 9968   | 9986   | 9996   | 9999   |
| 2   | 759   | 95     | 32     | 14     | 4      | 1      |
| 3   | 139   | 4      | 0      | 0      | 0      | 0      |
| 4   | 41    | 1      | 0      | 0      | 0      | 0      |
| 5–8 | 17    | 0      | 0      | 0      | 0      | 0      |

TABLE 4
*Empirical upper tail probabilities $P(S_{\hat{k}} \geq \chi^2_{5;\alpha})$ (upper line) and $P(S^*_{\hat{k}} \geq \chi^2_{5;\alpha})$ (lower line) from samples of size $n$ from the uniform distribution on $\mathbb{R}P^2$. The estimates are based on 10,000 simulations. The modified statistic $S^*_{\hat{k}}$ is defined in (5.2). Bold figures indicate values in the interval $\alpha \pm 2\sqrt{\alpha(1-\alpha)/10{,}000}$*

| $\alpha$ | $n=5$ | $n=10$ | $n=15$ | $n=20$ | $n=25$ | $n=30$ |
|---|---|---|---|---|---|---|
| 0.10 | 0.162 | **0.099** | **0.100** | **0.098** | **0.098** | **0.100** |
|      | 0.116 | **0.099** | **0.094** | **0.099** | **0.095** | **0.101** |
| 0.05 | 0.123 | **0.051** | **0.052** | **0.049** | **0.048** | **0.048** |
|      | 0.028 | 0.042 | 0.042 | **0.045** | 0.040 | **0.046** |
| 0.01 | 0.100 | 0.018 | 0.013 | **0.009** | **0.010** | **0.009** |
|      | 0.000 | **0.009** | 0.006 | 0.006 | 0.005 | **0.008** |

$$(5.3) \quad S^*_{\hat{k}} = \left\{1 + \frac{5.496 - 0.636 S_{\hat{k}} + 0.018 S^2_{\hat{k}}}{n}\right\} S_{\hat{k}} \qquad \text{for } M = SO(3)$$

Table 5
*Empirical distribution of $\hat{k}$ (based on 10,000 simulations) from samples of size n from the uniform distribution on $SO(3)$*

|     | $n=5$ | $n=10$ | $n=15$ | $n=20$ | $n=25$ | $n=30$ |
|-----|-------|--------|--------|--------|--------|--------|
| 1   | 9790  | 9997   | 10000  | 10000  | 10000  | 10000  |
| 2   | 189   | 3      | 0      | 0      | 0      | 0      |
| 3–4 | 21    | 0      | 0      | 0      | 0      | 0      |

Table 6
*Empirical upper tail probabilities $P(S_{\hat{k}} \geq \chi^2_{9;\alpha})$ (upper line) and $P(S^*_{\hat{k}} \geq \chi^2_{9;\alpha})$ (lower line) from samples of size n from the uniform distribution on $SO(3)$. The estimates are based on 10,000 simulations. The modified statistic $S^*_{\hat{k}}$ is defined in (5.3). Bold figures indicate values in the interval $\alpha \pm 2\sqrt{\alpha(1-\alpha)/10{,}000}$*

| $\alpha$ | $n=5$ | $n=10$ | $n=15$ | $n=20$ | $n=25$ | $n=30$ |
|----------|-------|--------|--------|--------|--------|--------|
| 0.10     | **0.094** | 0.087 | 0.093 | 0.093 | **0.097** | **0.095** |
|          | **0.095** | 0.089 | **0.098** | **0.096** | **0.097** | **0.095** |
| 0.05     | **0.055** | 0.042 | **0.045** | 0.044 | **0.047** | **0.047** |
|          | 0.043 | 0.037 | **0.045** | **0.046** | **0.045** | **0.045** |
| 0.01     | 0.027 | **0.008** | **0.010** | **0.008** | **0.009** | **0.009** |
|          | 0.024 | **0.008** | **0.009** | 0.012 | **0.009** | **0.009** |

(obtained by cubic regression of the empirical quantiles of the limiting distribution on the empirical quantiles of $S_{\hat{k}}$) is reasonable at the 10% and 5% levels for $n \geq 15$ for $S^2$, $\mathbb{R}P^2$ and $SO(3)$.

Table 7
*Empirical power $P(S \geq \chi^2_{3;\alpha})$ (based on 10,000 simulations) of tests of uniformity on $S^2$ against alternative (5.4) of mixture of Fisher distributions. $F_n$ denotes Giné's statistic; $S^*_{\hat{k}}$ is defined in (5.1)*

| $\alpha$ | $S$ | $n=5$ | $n=10$ | $n=15$ | $n=20$ | $n=25$ | $n=30$ |
|----------|-----|-------|--------|--------|--------|--------|--------|
| 0.10 | $F_n$ | 0.101 | 0.117 | 0.143 | 0.151 | 0.177 | 0.197 |
|      | $S_{\hat{k}}$ | 0.301 | 0.201 | 0.196 | 0.190 | 0.205 | 0.215 |
|      | $S^*_{\hat{k}}$ | 0.103 | 0.176 | 0.184 | 0.183 | 0.200 | 0.210 |
| 0.05 | $F_n$ | 0.050 | 0.060 | 0.070 | 0.078 | 0.096 | 0.106 |
|      | $S_{\hat{k}}$ | 0.261 | 0.156 | 0.145 | 0.147 | 0.162 | 0.171 |
|      | $S^*_{\hat{k}}$ | 0.000 | 0.127 | 0.137 | 0.138 | 0.157 | 0.165 |
| 0.01 | $F_n$ | 0.000 | 0.000 | 0.000 | 0.000 | 0.000 | 0.000 |
|      | $S_{\hat{k}}$ | 0.182 | 0.124 | 0.114 | 0.111 | 0.128 | 0.137 |
|      | $S^*_{\hat{k}}$ | 0.000 | 0.000 | 0.104 | 0.107 | 0.124 | 0.133 |

A common alternative to uniformity on $S^2$ is a mixture of two Fisher distributions with modes that are antipodal, so that the probability density function is

(5.4) $$f(\mathbf{x}; \boldsymbol{\mu}, \kappa) = \frac{\kappa}{2\sinh\kappa}(\exp\{\kappa\boldsymbol{\mu}'\mathbf{x}\} + \exp\{-\kappa\boldsymbol{\mu}'\mathbf{x}\})$$

with $\boldsymbol{\mu} \in S^2$. The main Sobolev test in common use on $S^2$ that is consistent against all alternatives is Giné's [7] $F_n$ test. (See also page 209 of [19].) Table 7 gives the empirical power (based on 10,000 simulations) of some tests of uniformity against the alternative of density (5.4) with $\kappa = 2$. The tests reject uniformity for $P(S \geq \chi^2_{3;\alpha})$, where $S$ is $F_n$, the statistic $S_{\hat{k}}$ of the data-driven test, and the modification $S^*_{\hat{k}}$ defined in (5.1). The data-driven test is definitely more powerful than Giné's test for the small sample sizes considered.

## APPENDIX: PROOF OF PROPOSITION 2.1

This proof is based on that of the asymptotic behavior of the heat kernel given in Chapter 5 of [23].

Let $D$ denote the de Rham operator $d + d^*$ acting on the complex exterior forms on $M$. Then $D^2 = \Delta$ is the Laplacian operator on the complex exterior forms on $M$. For $T = 1, 2, \ldots$, choose a smooth function $\varphi_T : \mathbb{R} \to \mathbb{R}$ such that

$$\varphi_T(s) = \begin{cases} 1, & \text{if } |s| \leq (2\lambda_T + \lambda_{T+1})/3, \\ 0, & \text{if } |s| \geq (\lambda_T + 2\lambda_{T+1})/3. \end{cases}$$

The kernel of $\varphi_T(\Delta)$ is the spectral function $e(\cdot, \cdot, T)$ given by (2.8). Given distinct points $x$ and $y$ in $M$, define $\delta$ as the Riemannian distance between $x$ and $y$. Choose a smooth function $\psi : \mathbb{R} \to \mathbb{R}$ such that

(A.1)
$$|t| \leq \delta/3 \implies \psi(t) = 1,$$
$$|t| \geq 2\delta/3 \implies \psi(t) = 0.$$

Since $\varphi_T$ and $\psi$ are in the class $\mathcal{S}(\mathbb{R})$ of rapidly decreasing functions with all derivatives rapidly decreasing, their Fourier transforms $\hat{\varphi}_T$ and $\hat{\psi}$ are in $\mathcal{S}(\hat{\mathbb{R}})$, and so there are functions $f_{1,T}$ and $f_{2,T}$ in $\mathcal{S}(\mathbb{R})$ with Fourier transforms $\hat{f}_{1,T}$ and $\hat{f}_{2,T}$ satisfying

(A.2)
$$\hat{f}_{1,T}(\lambda) = \frac{1}{\nu_T}\hat{\varphi}_T(\lambda)\psi(\lambda),$$
$$\hat{f}_{2,T}(\lambda) = \frac{1}{\nu_T}\hat{\varphi}_T(\lambda)(1 - \psi(\lambda)).$$

Then

$$f_{1,T}(\Delta) + f_{2,T}(\Delta) = \frac{1}{\nu_T}\varphi_T(\Delta).$$

Since the functions $f_{1,T}$ and $f_{2,T}$ are in the space $\mathcal{R}(\mathbb{R})$ of rapidly decreasing functions, the operators $f_{1,T}(\Delta)$ and $f_{2,T}(\Delta)$ have smooth kernels (see Proposition 5.8 of [23]). Denote these kernels by $e_{1,T}$ and $e_{2,T}$. As

$$f_{1,T}(\Delta) = \frac{1}{\nu_T} \int_{-\infty}^{\infty} \hat{\varphi}_T(\lambda^2) e^{i\lambda D} \psi(\lambda) \, d\lambda,$$

it follows from (A.1) and the fact that $D$ has unit propagation speed (cf. Proposition 5.5 of [23]) that $e_{1,T}$ is supported within $2\delta/3$ of the diagonal. Thus

$$\text{(A.3)} \qquad \frac{e(x,y,T)}{\nu_T} = e_{2,T}(x,y).$$

It follows from (A.2) that

$$\left| \hat{f}_{2,T}(\lambda) - \frac{2\sin(c_T \lambda)}{\lambda \nu_T}(1 - \psi(\lambda)) \right| \leq \frac{2(\lambda_{T+1} - \lambda_T)}{3\nu_T},$$

where $c_T = (2\lambda_T + \lambda_{T+1})/3$, and so $\hat{f}_{2,T}$ tends to 0 in $\mathcal{S}(\hat{\mathbb{R}})$ as $T \to \infty$. By Fourier theory, $f_{2,T} \to 0$ in $\mathcal{S}(\mathbb{R})$ and so in $\mathcal{R}(\mathbb{R})$. Thus (see Proposition 5.8 of [23]) $e_{2,T} \to 0$ as $T \to \infty$. By (A.3) this gives (2.14).

**Acknowledgment.** I am grateful to Jørgen Ellegaard Andersen for suggesting that Proposition 2.1 could be proved by adapting a proof of the asymptotic behavior of the heat kernel.

SCHOOL OF MATHEMATICS AND STATISTICS
UNIVERSITY OF ST. ANDREWS
NORTH HAUGH
ST. ANDREWS KY16 9SS
UNITED KINGDOM
E-MAIL: pej@st-andrews.ac.uk